\def\date{21.2.08}

\documentclass[12pt]{article}

\usepackage{amssymb}
\usepackage{amsmath}

\newcommand{\msk}{\vskip 5mm}  
  
\font\tengoth=eufm10 at 10pt
\font\sevengoth=eufm7 at 6pt
\newfam\gothfam
\textfont\gothfam=\tengoth
\scriptfont\gothfam=\sevengoth

\def\frak{\fam\gothfam\tengoth}

\newcommand{\g}{{\mathfrak g}}

\newcommand{\h}{{\mathfrak h}}

\newcommand{\fn}{{\mathfrak n}}

\renewcommand{\:}{\colon}
\newcommand{\0}{{\bf 0}}

\renewcommand{\phi}{\varphi}
\newcommand{\dd}{{\tt d}} 

\newcommand{\trile}{\trianglelefteq}
\newcommand{\subeq}{\subseteq}

\newcommand{\into}{\hookrightarrow}

\def\onto{\to\mskip-14mu\to} 

\newcommand{\1}{{\mathbf 1}}
\newcommand{\N}{{\mathbb N}}
\newcommand{\Z}{{\mathbb Z}}
\newcommand{\R}{{\mathbb R}}
\newcommand{\C}{{\mathbb C}}

\newcommand{\K}{{\mathbb K}}

\newcommand{\T}{{\mathbb T}}
\newcommand{\V}{{\mathbb V}}

\renewcommand{\hat}{\widehat}

\renewcommand{\L}{\mathop{\bf L{}}\nolimits}

\newcommand{\GL}{\mathop{{\rm GL}}\nolimits}
\newcommand{\PGL}{\mathop{{\rm PGL}}\nolimits}

\newcommand{\Gau}{\mathop{{\rm Gau}}\nolimits}

\newcommand{\aut} {\mathop{{\mathfrak{aut}}}\nolimits}

\newcommand{\gau} {\mathop{{\mathfrak{gau}}}\nolimits}
\newcommand{\gl}  {\mathop{{\mathfrak{gl} }}\nolimits}

\newcommand{\ad}{\mathop{{\rm ad}}\nolimits}
\newcommand{\Ad}{\mathop{{\rm Ad}}\nolimits}

\newcommand{\Hom}{\mathop{{\rm Hom}}\nolimits}

\newcommand{\Hol}{\mathop{{\rm Hol}}\nolimits}

\newcommand{\Aut}{\mathop{{\rm Aut}}\nolimits}

\newcommand{\Diff}{\mathop{{\rm Diff}}\nolimits}

\newcommand{\End}{\mathop{{\rm End}}\nolimits}
\newcommand{\id}{\mathop{{\rm id}}\nolimits}

\renewcommand{\dim}{\mathop{{\rm dim}}\nolimits}

\newcommand{\ev}{\mathop{{\rm ev}}\nolimits}
\newcommand{\evol}{\mathop{{\rm evol}}\nolimits}
\newcommand{\der}{\mathop{{\rm der}}\nolimits}

\newcommand{\nin}{\noindent} 
\newcommand{\oline}{\overline}

\newcommand{\res}{\vert}

\newcommand{\ssssarr}{\hbox to 15pt{\rightarrowfill}}
\newcommand{\sssarr}{\hbox to 20pt{\rightarrowfill}}
\newcommand{\ssarr}{\hbox to 30pt{\rightarrowfill}}
\newcommand{\sarr}{\hbox to 40pt{\rightarrowfill}}
\newcommand{\arr}{\hbox to 60pt{\rightarrowfill}}
\newcommand{\larr}{\hbox to 60pt{\leftarrowfill}}
\newcommand{\Arr}{\hbox to 80pt{\rightarrowfill}}
\newcommand{\mapdown}[1]{\Big\downarrow\rlap{$\vcenter{\hbox{$\scriptstyle#1$}}$}}

\newcommand{\sssmapright}[1]{\smash{\mathop{\sssarr}\limits^{#1}}}

\newcommand{\smapright}[1]{\smash{\mathop{\sarr}\limits^{#1}}}

\def\theoremname{Theorem}
\def\propositionname{Proposition}
\def\corollaryname{Corollary}
\def\lemmaname{Lemma}
\def\remarkname{Remark}
\def\conjecturename{Conjecture} 

\def\definitionname{Definition}
\def\exercisename{Exercise}
\def\examplename{Example}
\def\examplesname{Examples}
\def\problemname{Problem}
\def\problemsname{Problems}
\def\proofname{Proof}

\def\satzname{Satz} 
\def\koroname{Korollar}
\def\folgname{Folgerung}
\def\bemerkname{Bemerkung}
\def\aufgname{Aufgabe}

\def\beisname{Beispiel}
\def\beissname{Beispiele}
\def\bewname{Beweis}

\def\@thmcounter#1{\noexpand\arabic{#1}}
\def\@thmcountersep{}
\def\@begintheorem#1#2{\it \trivlist \item[\hskip 
\labelsep{\bf #1\ #2.\quad}]}
\def\@opargbegintheorem#1#2#3{\it \trivlist
      \item[\hskip \labelsep{\bf #1\ #2.\quad{\rm #3}}]}
\makeatother
\newtheorem{theo}{\theoremname}[section]
\newtheorem{propo}[theo]{\propositionname}
\newtheorem{coro}[theo]{\corollaryname}
\newtheorem{lemm}[theo]{\lemmaname}

\newenvironment{theorem}{\begin{theo}\it}{\end{theo}}

\newenvironment{proposition}{\begin{propo}\it}{\end{propo}}

\newenvironment{lemma}{\begin{lemm}\it}{\end{lemm}}

\newtheorem{rem}[theo]{\remarkname}

\newenvironment{remark}{\begin{rem}\rm}{\end{rem}}

\newtheorem{stepnow}[theo]{}

\newtheorem{defin}[theo]{\definitionname} 

\newenvironment{definition}{\begin{defin}\rm}{\end{defin}}

\newtheorem{exer}[theo]{\exercisename}

\newtheorem{ex}[theo]{\examplename}

\newenvironment{example}{\begin{ex}\rm}{\end{ex}}

\newtheorem{exs}[theo]{\examplesname}

\newenvironment{examples}{\begin{exs}\rm}{\end{exs}}

\newtheorem{conj}[theo]{\conjecturename}

\newtheorem{pr}[theo]{\problemname}

\newtheorem{prs}[theo]{\problemsname}

\newcommand{\qed}{{\unskip\nobreak\hfil\penalty50\hskip .001pt \hbox{}
          \nobreak\hfil
          \vrule height 1.2ex width 1.1ex depth -.1ex
           \parfillskip=0pt\finalhyphendemerits=0\medbreak}\rm}

{\hfill\qed\end{trivlist}}
\newenvironment{proof}{\begin{trivlist}\item[\hskip%
\labelsep{\bf\proofname.\quad}]}%
{\hfill\qed\end{trivlist}}

\newenvironment{Proof*}{\begin{trivlist}\item[\hskip%
\labelsep{\bf\proofname.\quad}]}%
{\end{trivlist}}

\newcommand{\pmat}[1]{\begin{pmatrix} #1 \end{pmatrix}}

{\hfill\qed\end{trivlist}}

\newenvironment{beweis*}{\begin{trivlist}\item[\hskip%
\labelsep{\bf\bewname.\quad}]}%
{\end{trivlist}}

\newtheorem{satzn}[theo]{\satzname}

\newtheorem{koro}[theo]{\koroname}

\newtheorem{folg}[theo]{\folgname}

\newtheorem{bem}[theo]{\bemerkname}

\newtheorem{aufg}[theo]{\aufgname}

\newtheorem{aufgn}[theo]{\aufgname}

\newtheorem{beis}[theo]{\beisname}

\newtheorem{beiss}[theo]{\beissname}

\begin{document} 


\def\V{{{\mathbb V}}} 
\def\GaL{\mathop{\rm \Gamma L}\nolimits}
\def\Car{\mathop{\rm Car}\nolimits}
\def\Emb{\mathop{\rm Emb}\nolimits}
\def\Fr{\mathop{\rm Fr}\nolimits}
\def\gal{\mathop{\rm \gamma {\frak l}}\nolimits}
\def\dend{\mathop{\rm DEnd}\nolimits}
\def\Idem{\mathop{\rm Idem}\nolimits}
\def\W{{\mathbb W}}
\def\dd{{\tt d}} 
\newcommand{\n}{\fn}
\newcommand{\Inv}{\mathop{\rm Inv}\nolimits} 
\newcommand{\HC}{\mathop{\rm HC}\nolimits} 

\title{Lie group extensions associated to projective modules 
of continuous inverse algebras} 
\author{K.-H. Neeb} 

\maketitle

\begin{abstract} 
{We call a unital locally convex algebra $A$ a continuous inverse 
algebra if its unit group $A^\times$ is open and inversion is a continuous 
map. For any smooth action of a, 
possibly infinite-dimensional, connected 
Lie group $G$ on a continuous inverse algebra $A$ by automorphisms 
and any finitely generated projective right $A$-module $E$, we construct 
a Lie group extension $\hat G$ of $G$ by the group 
$\GL_A(E)$ of automorphisms of the $A$-module $E$. This Lie group 
extension is a ``non-commutative'' version of the group 
$\Aut(\V)$ of automorphism of a vector bundle over a compact manifold 
$M$, which arises for $G = \Diff(M)$, $A = C^\infty(M,\C)$ and 
$E = \Gamma\V$. We also identify the Lie algebra $\hat\g$ of $\hat G$ 
and explain how it is related to connections of the $A$-module $E$. \\
{\sl AMC Classification:} 22E65, 58B34}
\end{abstract}


\section*{Introduction} 

In \cite{ACM89} it is shown that for a finite-dimensional 
$K$-principal bundle 
$P$ over a compact manifold $M$, the group $\Aut(P)$ of all bundle automorphisms 
carries a natural Lie group structure whose Lie algebra is the Fr\'echet--Lie algebra 
of ${\cal V}(P)^K$ of $K$-invariant smooth vector fields on $M$. 
This applies in particular to the group $\Aut(\V)$ of automorphisms of a 
finite-dimensional vector bundle with fiber $V$ 
because this group can be identified with the 
automorphisms group of the corresponding frame bundle $P = \Fr \V$ which is a 
$\GL(V)$-principal bundle. 

In this paper, we turn to variants of the Lie groups 
$\Aut(\V)$ arising in non-commutative geometry. 
In view of \cite{Ko76}, the group $\Aut(\V)$ can be identified with the 
group of semilinear automorphisms of the $C^\infty(M,\R)$-module $\Gamma(\V)$ of 
smooth sections of $\V$, which, according to Swan's Theorem, is a finitely generated 
projective module. Here the gauge group $\Gau(\V)$ corresponds to the 
group of $C^\infty(M,\R)$-linear module isomorphisms. 

This suggests the following setup: Consider a unital locally convex algebra 
$A$ and a finitely generated projective right $A$-module $E$. When can we turn 
groups of semilinear automorphisms of $E$ into Lie groups? First of all, we have 
to restrict our attention to a natural class of algebras whose unit groups 
$A^\times$ carry natural Lie group structures, which is the case if 
$A^\times$ is an open subset of $A$ and the inversion map is continuous. 
Such algebras are called {\it continuous inverse algebras}, CIAs, for short.  
The Fr\'echet algebra $C^\infty(M,\R)$ is a CIA if and only if $M$ is compact. 
Then its automorphism group  $\Aut(C^\infty(M,\R)) \cong 
\Diff(M)$ carries a natural Lie group 
structure with Lie algebra ${\cal V}(M)$, the Lie algebra of smooth 
vector fields on $M$. 
Another important class of CIAs whose automorphism groups 
are Lie groups are 
smooth $2$-dimensional quantum tori with generic diophantine properties 
(cf.\ \cite{El86}, \cite{BEGJ89}). 
Unfortunately, in general, automorphism groups of 
CIAs do not always carry a natural Lie group structure, so that it is much more 
natural to consider triples $(A,G,\mu_A)$, where $A$ is a CIA, $G$ a 
possibly infinite-dimensional Lie group, and $\mu_A \: G \to \Aut(A)$ 
a group homomorphism defining a smooth action of $G$ on $A$. 

For any such triple $(A,G,\mu_A)$ and any finitely generated 
projective $A$-module $E$, the subgroup $G_E$ of all elements of $G$ for which 
$\mu_A(g)$ lifts to a semilinear automorphism of $E$ is an open subgroup. 
One of our main results (Theorem~\ref{thm:3.3}) is that we thus obtain a Lie group 
extension 
$$ \1 \to \GL_A(E) = \Aut_A(E) \to \hat G_E \to G_E \to \1, $$
where $\hat G_E$ is a Lie group acting smoothly on $E$ by semilinear automorphisms. 
For the special case where $M$ is a compact manifold, 
$A = C^\infty(M,\R)$, $E = \Gamma(\V)$ for a smooth vector bundle $\V$,  
and  $G = \Diff(M)$, the Lie group 
$\hat G$ is isomorphic to 
the group $\Aut(\V)$ of automorphisms of the vector bundle $\V$, but our 
construction contains a variety of other interesting settings. 
From a different perspective, the Lie group structure on $\hat G$ also 
tells us about possible smooth actions of Lie groups on finitely generated 
projective $A$-modules by semilinear maps which are compatible with a 
smooth action on the algebra $A$. 

A starting point of our construction is the observation that the connected components 
of the set $\Idem(A)$ of idempotents of a CIA coincide with the orbits 
of the identity component $A^\times_0$ of $A^\times$ under the conjugation 
action. Using the natural manifold structure on $\Idem(A)$ 
(cf.\ \cite{Gram84}), the action of $A^\times$ on $\Idem(A)$ even is a 
smooth Lie group action. 

On the Lie algebra side, the semilinear automorphisms of $E$ correspond
to the Lie algebra $\dend(E)$ of derivative endomorphisms, i.e., those 
endomorphisms $\phi \in \End_\K(E)$ for which there is a continuous 
derivation $D \in \der(A)$ 
with $\phi(s.a) = \phi(s).a + s.(D.a)$ for $s \in E$ and $a \in A$. 
The set $\hat\dend(E)$ of all pairs $(\phi, D) \in \End_\K(E)\times \der(A)$ 
satisfying this condition is a Lie algebra and we obtain a Lie algebra 
extension 
$$ \0 \to \End_A(E) = \gl_A(E) \into \hat\dend(E) \onto \der(A) \to \0. $$
Pulling this extension  back via the Lie algebra homomorphism 
$\g\to \der(A)$ induced by the action of $G$ on $A$ yields the 
Lie algebra $\hat\g$ of the group $\hat G$ from above 
(Proposition~\ref{prop:4.7}). 

In Section~\ref{sec:5} we briefly discuss the relation between 
linear splittings of the Lie algebra extension $\hat \g$ and covariant 
derivatives in the context of non-commutative geometry 
(cf.\ \cite{Co94}, \cite{MMM95}, \cite{DKM90}). 


{\bf Thanks:} We thank Hendrik Grundling for reading erlier versions 
of this paper and for numerous remarks which lead to several improvements 
of the presentation.

\subsection*{Preliminaries and notation} 

Throughout this paper we write $I := [0,1]$ for the unit interval in $\R$ and 
$\K$ either denotes $\R$ or $\C$. 
A locally convex space $E$ is said to be {\it Mackey complete} if 
each smooth curve $\gamma \: I \to E$ has a (weak) integral in $E$.
For a more detailed discussion of Mackey completeness and 
equivalent conditions, we refer to \cite[Th.~2.14]{KM97}. 
 
A {\it Lie group} $G$ is a group equipped with a 
smooth manifold structure modeled on a locally convex space 
for which the group multiplication and the 
inversion are smooth maps. We write $\1 \in G$ for the identity element and 
$\lambda_g(x) = gx$, resp., $\rho_g(x) = xg$ for the left, resp.,
right multiplication on $G$. Then each $x \in T_\1(G)$ corresponds to
a unique left invariant vector field $x_l$ with 
$x_l(g) := d\lambda_g(\1).x, g \in G.$
The space of left invariant vector fields is closed under the Lie
bracket of vector fields, hence inherits a Lie algebra structure. 
In this sense we obtain on $\g := T_\1(G)$ a continuous Lie bracket which
is uniquely determined by $[x,y]_l = [x_l, y_l]$ for $x,y \in \g$. 
We shall also use the functorial 
notation $\L(G)$ for the Lie algebra of $G$ and, accordingly, 
$\L(\phi) = T_\1(\phi)\: \L(G_1) \to \L(G_2)$ 
for the Lie algebra morphism associated to 
a morphism $\phi \: G_1 \to G_2$ of Lie groups. 

A Lie group $G$ is called {\it regular} if for each 
$\xi \in C^\infty(I,\g)$, the initial value problem 
$$ \gamma(0) = \1, \quad \gamma'(t) = \gamma(t).\xi(t) = T(\lambda_{\gamma(t)})\xi(t) $$
has a solution $\gamma_\xi \in C^\infty(I,G)$, and the
evolution map 
$$ \evol_G \: C^\infty(I,\g) \to G, \quad \xi \mapsto \gamma_\xi(1) $$
is smooth (cf.\ \cite{Mil84}). 
For a locally convex space $E$, 
the regularity of the Lie group $(E,+)$ is equivalent 
to the Mackey completeness of $E$ (\cite[Prop.~II.5.6]{Ne06}). 
We also recall that for each regular 
Lie group $G$ its Lie algebra $\g$ is Mackey complete 
and that all Banach--Lie groups are regular (cf.\ \cite[Rem.~II.5.3]{Ne06} 
and \cite{GN08}). 

A smooth map $\exp_G \: \L(G) \to G$  is called an {\it exponential function} 
if each curve $\gamma_x(t) := \exp_G(tx)$ is a one-parameter group 
with $\gamma_x'(0)= x$. The Lie group $G$ is said to be {\it locally exponential} 
if it has an exponential function for which there is an open $0$-neighborhood 
$U$ in $\L(G)$ mapped diffeomorphically by $\exp_G$ onto an open subset of $G$. 

If $A$ is an associative algebra with unit, we write $\1$ for the identity element, 
$A^\times$ for its group of units, $\Idem(A) = \{ p \in A \: p^2 = p\}$ for its 
set of idempotents and $\eta_A(a)= a^{-1}$ for the inversion map $A^\times \to A$. 
A homomorphism $\rho \: A \to B$ is unital algebras is called {\it isospectral} 
if $\rho^{-1}(B^\times) = A^\times$. We write $\GL_n(A) := M_n(A)^\times$ 
for the unit group of the unital algebra $M_n(A)$ of $n \times n$-matrices 
with entries in $A$. 

Throughout, $G$ denotes a (possibly infinite-dimensional) Lie 
group , $A$ a Mackey complete unital continuous inverse algebra (CIA for short) and 
$G \times A \to A, (g,a) \mapsto g.a = \mu_A(g)(a)$ is a smooth action 
of $G$ on~$A$. 

\section{Idempotents and finitely generated projective modules} \label{sec:1}

The set $\Idem(A)$ of idempotents of a CIA $A$ plays a central 
role in (topological) $K$-theory. In \cite[Satz~2.13]{Gram84}, 
Gramsch shows that this set always 
carries a natural manifold structure, 
which implies in particular that its connected 
components are open subsets. 
Since we shall need it in the following, we briefly recall some basic facts 
on $\Idem(A)$ (cf.\ \cite{Gram84}; see also \cite[Sect.~4]{Bl98}).  

\begin{proposition} \label{prop:1.1}
For each $p \in \Idem(A)$, the set 
$$ U_p := \{ q \in \Idem(A) \: pq + (\1-p)(\1-q) \in A^\times \} $$
is an open neighborhood of $p$ in $\Idem(A)$ 
and, for each $q \in U_p$, the element 
$$ s_q := pq + (\1-p)(\1-q) \in A^\times 
\quad \hbox{ satisfies } \quad s_q q s_q^{-1} = p. $$
The connected component of $p$ in $\Idem(A)$ coincides with the 
orbit of the identity component $A^\times_0$ of $A^\times$ 
under the conjugation action $A^\times \times \Idem(A) \to \Idem(A), 
(g,p) \mapsto c_g(p) := g.p := gpg^{-1}$. 
\end{proposition}

\begin{proof}
Since the map $q \mapsto pq + (\1-p)(\1-q)$ is continuous, 
it maps $p$ to $\1$ and since $A^\times$ is open, $U_p$ is an open neighborhood 
of $p$. Hence, for each $q \in U_p$, the element $s_q$ is invertible, 
and a trivial calculation shows that $s_q q = p s_q$. 

If $q$ is sufficiently close to $p$, then $s_q \in A^\times_0$ because 
$s_p = \1$ and $A^\times_0$ is an open neighborhood of $\1$ in $A$ 
(recall that $A$ is locally convex and $A^\times$ is open). 
This implies that $q = s_q^{-1}ps_q$ lies in the orbit 
$\{ c_g(p) \: g \in A^\times_0\}$ of $p$ under $A^\times_0$. 
Conversely, since the orbit map $A^\times \to \Idem(A), 
g \mapsto c_g(p)$ is continuous, 
it maps the identity component  $A_0^\times$ into the connected component 
of $p$ in $\Idem(A)$. 
\end{proof}

\begin{lemma} \label{lem:1.2}
Let $A$ be a CIA, $n \in \N$, and $p = p^2 \in M_n(A)$ an idempotent. 
Then the following assertions hold: 
\begin{description}
\item[\rm(1)] The subalgebra $p M_n(A)p$ is a CIA with identity element $p$. 
\item[\rm(2)] The unit group $(p M_n(A)p)^\times$ is a Lie group. 
\end{description}
\end{lemma}

\begin{proof}
Since $M_n(A)$ also is a CIA (\cite{Swa62}; \cite{Gl02}) 
which is Mackey complete if $A$ is so, it suffices to prove the assertion for $n = 1$. 

(1) From the decomposition of the identity $\1$ as a sum $\1 = p + (\1- p)$ 
of two orthogonal idempotents, we obtain the direct sum decomposition 
$$ A = pAp \oplus p A (\1-p) \oplus (\1- p) A p \oplus (\1-p)A(\1-p). $$

We claim that an element $a \in pAp$ is invertible in this algebra if and only 
if the element $a + (\1 - p)$ is invertible in $A$. 
In fact, if $b \in pAp$ satisfies $ab = ba = p$, then 
$$ (a + (\1-p))(b + (\1-p)) = ab + (\1- p) = \1 
=  (b + (\1-p))(a + (\1-p)). $$
If, conversely, $c \in A$ is an inverse of $a + (\1-p)$ in $A$, then 
$ca + c(\1-p) = \1 = ac + (\1-p)c$
leads after multiplication with $p$ to 
$ca = p = ac$, which implies 
$(pcp)a = p  = a(pcp)$, so that $pcp$ is an inverse of $a$ in $pAp$. 
The preceding argument implies in particular that 
$(pAp)^\times = pAp \cap (A^\times - (\1 - p))$ is an open subset of $pAp$, 
and that the inversion map 
$$ \eta_{pAp} \: (pAp)^\times  \to pAp, \quad 
a \mapsto a^{-1} = \eta_A(a + \1 - p) - (\1 -p) $$
is continuous. 

(2) is an immediate consequence of (1) (cf.\ \cite{Gl02}, 
\cite{Ne06}, Ex.~II.1.4, Th.~IV.1.11). 
\end{proof}

Let $E$ be a finitely generated projective right $A$-module. 
Then there is some $n \in \N$ and an idempotent $p = p^2 \in M_n(A)$ with 
$E \cong p A^n$, where $A$ acts by multiplication on the right. Conversely, 
for each idempotent $p \in \Idem(M_n(A))$, the right submodule $pA^n$ of $A^n$ is 
finitely generated (as a quotient of $A^n$) and projective because it is a direct 
summand of the free module $A^n \cong p A^n \oplus (\1-p)A^n$. 
The following lemma provides some information on 
$A$-linear maps between such modules. 

\begin{lemma} \label{lem:1.3}
Let $p,q\in \Idem(M_n(A))$ be two idempotents. 
Then the following assertions hold: 
\begin{description}
\item[(1)] The map 
$x \mapsto \lambda_x\res_{pA^n}$ (left multiplication) yields a bijection 
$$ \alpha_{p,q} \: 
q M_n(A) p = \{ x \in M_n(A) \: qx = x, xp = x\} \to \Hom_A(pA^n,qA^n). $$
\item[(2)] $pA^n \cong qA^n$ if and only if there are $x,y \in M_n(A)$ with 
$xy = q$ and $yx = p$. If, in particular, $q = gpg^{-1}$ for some 
$g \in \GL_n(A)$, then $x := gp$ and $y := g^{-1}$ satisfy 
$xy = q$ and $yx = p$. 
\item[(3)] $pA^n \cong qA^n$ if and only if there are $x \in qM_n(A)p$ and 
$y \in pM_n(A)q$ with $xy = q$ and $yx = p$. 
\item[(4)] If $pA^n \cong qA^n$, then there exists an element 
$g \in \GL_{2n}(A)$ with $g \tilde p g^{-1}= \tilde q$, where 
$$ \tilde p =\pmat{p & 0 \cr 0 & 0 \cr} 
\quad \hbox{ and } \quad 
\tilde q =\pmat{q & 0 \cr 0 & 0 \cr}. $$
\end{description}
\end{lemma}

\begin{proof}
(1) Each element of $\End_A(A^n)$ is given by left multiplication 
with a matrix, so that $M_n(A) \cong \End_A(A^n)$. Since $pA^n$ and $qA^n$ 
are direct summands of $A^n$, each element of $\Hom(pA^n,qA^n)$ can be realized 
by left multiplication with a matrix, and we have the direct sum decomposition 
\begin{align*}
\End_A(A^n) &\cong \Hom(pA^n, qA^n) \oplus  \Hom(pA^n, (\1-q)A^n) \cr
&\oplus \Hom((\1-p)A^n, qA^n) 
\oplus \Hom((\1-p)A^n, (\1-q)A^n),
\end{align*}
which corresponds to the direct sum decomposition 
$$ M_n(A) \cong 
q M_n(A) p \oplus (\1-q)M_n(A) p \oplus 
q M_n(A) (\1-p) \oplus (\1-q)M_n(A) (\1-p). $$
Now the assertion follows from 
$qM_n(A)p = \{ x \in M_n(A) \: qx=x, xp= x\}.$

(2), (3) If 
$pA^n$ and $qA^n$ are isomorphic, there is some $x \in qM_n(A)p \cong \Hom(pA^n,qA^n)$ 
for which $\lambda_x \: pA^n \to qA^n, s \mapsto xs$ is an isomorphism. 
Writing $\lambda_x^{-1}$ as $\lambda_y$ for some $y \in \Hom(qA^n, pA^n) \cong pM_n(A)q$, 
we get 
$$ p = \lambda_y \circ \lambda_x(p) = yxp = yx 
\quad \hbox{ and } \quad 
q = \lambda_x \circ \lambda_y(q) = xyq = xy. $$

If, conversely, $p = yx$ and $q = xy$ hold for some $x,y \in M_n(A)$, 
then $p^2 = p$ implies 
$p = yx = yxyx = yqx$ and likewise $q = xy = xpy$, 
which leads to 
$$ (pyq)(qxp) = pyqxp=p^3 = p \quad \hbox{ and } \quad 
(qxp)(pyq) = qxpyq = q^3 = q. $$
Therefore, we also have $p = y'x'$ and $q=x'y'$ with 
$x' := qxp \in q M_n(A)p$ and $y' := pyq \in p M_n(A)q$. 
Then 
$\lambda_{x'} \: pA^n \to qA^n$ and 
$\lambda_{y'} \: qA^n \to pA^n$ are module homomorphisms with 
$\lambda_{x'} \circ \lambda_{y'} = \lambda_{x'y'} = \lambda_q = \id_{qA^n}$ and 
$\lambda_{y'} \circ \lambda_{x'} = \lambda_{y'x'} = \lambda_p = \id_{pA^n}$. 

(4) (cf.\ \cite[Prop.~4.3.1]{Bl98}) Pick $x,y$ as in (3). Let 
$$ \alpha := \pmat{ \1 - q & x \cr y & \1 -p \cr} \quad 
\hbox{ and } \quad \beta := \pmat{ \1 - p & p \cr p & \1 -p \cr} \in M_{2n}(A). $$
Then a direct calculation yields $\alpha^2 = \1 = \beta^2$. Therefore 
$z := \beta\alpha \in \GL_{2n}(A)$. Moreover, we have 
$$\alpha \tilde q \alpha^{-1} =  \pmat{ 0 & 0 \cr 0 & p\cr} 
\quad \hbox{ and } \quad 
\beta \pmat{ 0 & 0 \cr 0 & p\cr} \beta^{-1} = \tilde p, $$
so that $z \tilde q z^{-1} = \tilde p$. 
\end{proof}

\begin{proposition}
  \label{prop:1.4}
For each finitely generated projective right $A$-module $E$, 
we pick some idempotent $p \in M_n(A)$ with $E \cong p A^n$. Then 
we topologize $\End_A(E)$ by declaring 
$$\alpha_{p,p} \: pM_n(A)p \to \End_A(E), \quad 
x \mapsto \lambda_x\res_{pA^n} $$
to be a topological isomorphism. Then the 
algebra $\End_A(E)$ is a CIA and $\GL_A(E)$ is a Lie group. 
This topology does not depend on the choice of $p$ and 
if $A$ is Mackey complete, then $\GL_A(E)$ is locally exponential. 
\end{proposition}

\begin{proof}
We simply combine Lemma~\ref{lem:1.2} with Lemma~\ref{lem:1.3}(1) 
to see that 
we obtain a CIA structure on $\End_A(E)$, so that 
$\GL_A(E)$ is a Lie group which is locally exponential if $A$ is 
Mackey complete (\cite{Gl02}). 

To verify the independence of the topology on $\End_A(E)$ from $p$, 
we first note that for any matrix 
$$ \tilde p =\pmat{p & 0 \cr 0 & 0 \cr} \in M_N(A), \quad N > n, $$
we have a natural isomorphism 
$\tilde p M_N(A) \tilde p \cong p M_n(A) p,$
because all non-zero entries of matrices of the form 
$\tilde p X \tilde p$, $X \in M_N(A)$, 
lie in the upper left $(p \times p)$-submatrix and depend only on 
the corresponding entries of $X$. 

If $q \in \Idem(M_\ell(A))$ is another idempotent with 
$q A^\ell \cong E$, then the preceding argument shows that, after passing to $\max(n,\ell)$, 
we may w.l.o.g.\ assume that $\ell = n$. 
Then Lemma~\ref{lem:1.3} yields a $g \in \GL_{2n}(A)$ with 
$g\oline p g^{-1} = \oline q$, and then conjugation with $g$ 
induces a topological isomorphism 
$$ p M_n(A) p \cong 
\oline p M_{2n}(A) \oline p \smapright{c_g} 
\oline q M_{2n}(A) \oline q \cong q M_n(A) q. $$\\[-12mm]
\end{proof}

\begin{example} \label{ex:1.5}
 (a) Let $M$ be a smooth paracompact finite-dimensional 
manifold. We 
endow $A := C^\infty(M,\K)$ 
with the smooth compact open topology, i.e., the 
topology inherited by the natural embedding 
$$ C^\infty(M,\K)\into \prod_{k = 0}^\infty C(T^k M,T^k \K), \quad 
f \mapsto (T^k(f))_{k \in \N_0}, $$
where all spaces $C(T^k M,T^k \K)$ carry the compact open topology which 
coincides with the topology of uniform convergence on compact subsets. 

If $E$ is the space of smooth sections of a 
smooth vector bundle \break $q \: \V \to M$, then 
$E$ is a finitely generated projective $A$-module (\cite{Swa62}). 
The algebra $\End_A(E)$ is the space of smooth sections of the vector bundle
$\End(\V)$ and its unit group $\GL_A(E) \cong \Gau(\V)$ is the corresponding gauge group. 
We shall return to this class of examples below. 

(b) We obtain a similar picture if $A$ is the Banach algebra $C(X,\K)$, 
where $X$ is a compact space and 
$E$ is the space of continuous sections of a finite-dimensional topological 
vector bundle over $X$. Then $\End_A(E)$ is a Banach algebra, so that its unit group 
$\GL_A(E)$ is a Banach--Lie group. 
\end{example}

\section{Semilinear automorphisms of 
finitely generated projective modules} \label{sec:2}

In this section we take a closer look at the group 
$\GaL(E)$ of semilinear automorphism of a right $A$-module $E$. 
One of our main results, proved in Section~\ref{sec:3} below, asserts that 
if $E$ is a finitely generated projective module of a CIA $A$, 
certain pull-backs of this group by a smooth action of a Lie group $G$ 
on $A$ lead to a Lie group extension $\hat G$ of $G$ by the Lie group 
$\GL_A(E)$ 
(cf.\ Proposition~\ref{prop:1.4}) acting smoothly on $E$. 

The discussion in this section is very much inspired by Y.~Kosmann's paper 
\cite{Ko76}. 

\begin{definition}
  \label{def:2.1}
Let $E$ be a topological right module of the CIA $A$, i.e., 
we assume that the module structure $E \times A \to E, (s,a) \mapsto s.a =:\rho_E(a)s$ 
is a continuous bilinear map. 
We write $\End_A(E)$ for the algebra of continuous module 
endomorphisms of $E$ and 
$\GL_A(E) := \End_A(E)^\times$ for its group of units, the module 
automorphism group of $E$. For $A = \K$ we have in particular 
$\GL(E) = \GL_\K(E)$. 
The group $\GL_A(E)$ is contained in the group 
\begin{align*}
&\ \ \ \GaL(E) \\
&:= \{ \phi \in \GL_\K(E) \: (\exists \phi_A \in \Aut(A))
(\forall s \in E)(\forall a \in A)\ \phi(s.a) = \phi(s).\phi_A(a)\}\cr
&= \{ \phi \in \GL_\K(E) \: (\exists \phi_A \in \Aut(A))(\forall a \in A) 
\ \phi \circ \rho_E(a) =\rho_E(\phi_A(a)) \circ \phi\}. 
\end{align*}
of {\it semilinear automorphisms of $E$}, 
where we write $\Aut(A)$ for the group of topological automorphisms of~$A$.  
We put  
$$ \hat\GaL(E) := \{ (\phi,\phi_A) \in \GL_\K(E) \times \Aut(A) \: 
(\forall a \in A)\ \ \phi \circ \rho_E(a) =\rho_E(\phi_A(a)) \circ \phi\},  $$
where the multiplication is componentwise multiplication in the product group. 
In \cite{Har76}, the elements of $\hat\GaL(E)$ are called semilinear 
automorphisms and, for $A$ commutative, certain characteristic 
cohomology classes of $E$ are constructed for this group with values 
in differential forms over $A$. 
If the representation of $A$ on $E$ is faithful, then 
$\phi_A$ is uniquely determined by $\phi$, so that 
$\hat\GaL(E) \cong \GaL(E)$. 

The map 
$(\phi,\phi_A) \mapsto \phi_A$ defines a short exact sequence of groups 
$$ \1 \to \GL_A(E) \to \hat\GaL(E) \to \Aut(A)_E \to \1, $$
where $\Aut(A)_E$ denotes the image of the group 
$\hat\GaL(E)$ in $\Aut(A)$. 
\end{definition}

\begin{remark} \label{rem:2.2}
(a) For each $\psi \in \Aut(A)$, we define the 
corresponding {\it twisted module} 
$E^\psi$ by endowing the vector space $E$ with the new $A$-module 
structure defined by $s *_\psi a := s.\psi(a)$, i.e., 
$\rho_{E^\psi} = \rho_E \circ \psi$. Then a continuous linear map 
$\phi \: E \to E^\psi$ is a morphism of $A$-modules if and only if 
$\phi(s.a) = \phi(s).\psi(a)$ holds for all $s \in E$ and $a \in A$, i.e., 
$$ \phi \circ \rho_E(a) = \rho_{E^\psi}(a) \circ \phi \quad \hbox{ for } 
\quad a \in A. $$
Therefore $(\phi,\psi) \in \hat\GaL(E)$ is equivalent to 
$\phi \: E \to E^\psi$ being a module isomorphism. This shows that 
$$ \Aut(A)_E = \{ \psi \in \Aut(A) \: E^\psi \cong E\}.  $$

(b) Let $\psi \in \Diff(M)$ and $q \: \V \to M$ be a smooth vector bundle over $M$. 
We consider the pull-back vector bundle 
$$ \V^\psi := \psi^*\V := \{(x,v) \in M \times \V \: \psi(x) = q_\V(v)\} $$
with the bundle projection 
$q_\V^\psi \: \V^\psi \to M, (x,v) \mapsto x.$ 

If $s \: M \to \V^\psi$ is a smooth section, then 
$s(x) = (x, s'(\psi(x)))$, where $s' \: M \to \V$ is a smooth section, and this 
leads to an identification of the spaces of smooth sections of 
$\V$ and $\V^\psi$. 
For a smooth function $f \: M \to \K$ and $s \in \Gamma(\V^\psi)$, we have 
$(s.f)(x) = f(x)s(x)= (x, f(x) s'(\psi(x))),$
so that the corresponding right module structure on $E = \Gamma \V$ 
is given by $s' * f = s'.(\psi.f)$, 
where $\psi.f = f \circ \psi^{-1}$. This shows that 
$E^\psi = (\Gamma\V)^\psi \cong \Gamma(\V^\psi),$
i.e., the sections of the pull-back bundle form a twisted module. 

(c) Let $E$ be a finitely generated projective right $A$-module 
and $p \in \Idem(M_n(A))$ with $E \cong pA^n$. 
For $\psi \in \Aut(A)$ we write 
$\psi^{(n)}$ for the corresponding automorphisms 
of $A^n$, resp., $M_n(\psi)$ for the corresponding automorphism of 
$M_n(A)$. Then the map 
$\psi_{n} \: M_n(\psi)^{-1}(p) A^n \to p A^n$ 
induces a module isomorphism $M_n(\psi)^{-1}(p)A^n \cong (pA^n)^\psi$. 

(d) Let $\rho_E \: A \to \End(E)$ denote the representation of $A$ on $E$ defining 
the right module structure. Then, for each $a \in A^\times$, we have 
$(\rho_E(a), c_a^{-1}) \in \hat\GaL(E)$ because 
$\rho_E(a)(s.b) = s.ba = (s.a)(a^{-1}ba).$
\end{remark}

\begin{definition}
  \label{def:2.3}
Let $G$ be a group acting by automorphism on the group $N$ via $\alpha \: G \to \Aut(N)$. 
We call a map $f \: G \to N$ a {\it crossed homomorphism} if 
$$ f(g_1 g_2) = f(g_1) \alpha(g_1)(f(g_2)) 
\quad \hbox{ for } \quad g_1, g_2 \in G. $$
Note that $f$ is a crossed homomorphism if and only if 
$(f,\id_G) \: G \to N \rtimes_\alpha G$ is a group homomorphism. 
The set of all crossed homomorphisms $G \to N$ is denoted by $Z^1(G,N)$. 
The group $N$ acts naturally on $Z^1(G,N)$ by 
$$ (n * f)(g) := n f(g) \alpha(g)(n)^{-1} $$
and the set of $N$-orbits in $Z^1(G,N)$ is denoted $H^1(G,N)$. 
If $N$ is not abelian, $Z^1(G,N)$ and $H^1(G,N)$ do not carry a 
group structure; only the 
constant map $\1$ is a distinguished element of $Z^1(G,N)$, and 
its class $[\1]$ is distinguished in $H^1(G,N)$. The crossed homomorphisms 
in the class $[\1]$ are called {\it trivial}. They are of the form 
$f(g) = n\alpha(g)(n)^{-1}$ for some $n \in N$. 
\end{definition}

\begin{proposition}
  \label{prop:2.4}
Let $\rho_E \: A \to \End_\K(E)$ denote the action of $A$ 
on the right $A$-module $E$ and consider the action of the unit group $A^\times$ on the group $\GaL(E)$ by $\tilde\rho_E(a)(\phi) 
:= \rho_E(a)^{-1}\phi \rho_E(a).$ 
To each $\psi \in \Aut(A)$ we associate the function 
$$ C(\psi) \: A^\times \to \GaL(E), \quad 
a \mapsto \rho_E(\psi(a)a^{-1})^{-1} = \rho_E(\psi(a))^{-1} \rho_E(a). $$
Then $C(\psi) \in Z^1(A^\times, \GaL(E))$, and we thus get an exact sequence of 
pointed sets 
$$ \1 \to \GL_A(E) \to \hat\GaL(E) \to \Aut(A) \sssmapright{\oline C} 
 H^1(A^\times, \GaL(E)),$$
characterizing the subgroup $\Aut(A)_E$ as $\oline C^{-1}([\1])$. 
\end{proposition}

\begin{proof}
That $C(\psi)$ is a crossed homomorphism follows from 
\begin{align*}
C(\psi)(ab) 
&= \rho_E(\psi(ab))^{-1} \rho_E(ab)
= \rho_E(\psi(a))^{-1} \rho_E(\psi(b))^{-1} \rho_E(b)\rho_E(a) \cr
&= C(\psi)(a) \rho_E(a)^{-1} C(\psi)(b)\rho_E(a) 
= C(\psi)(a) \tilde\rho_E(a)\big(C(\psi)(b)\big).
\end{align*}
That the crossed homomorphism $C(\psi)$ is trivial means that there is a 
$\phi \in \GaL(E)$ with 
$$ C(\psi)(a) = \rho_E(\psi(a))^{-1} \rho_E(a) 
= \phi \rho_E(a)^{-1} \phi^{-1} \rho_E(a), $$
which means that 
$\rho_E(\psi(a)) \phi = \phi \rho_E(a)$ for $a \in A^\times.$
Since each CIA $A$ is generated by it unit group, which is an open subset, 
the latter relation is equivalent to $(\phi,\psi) \in \hat\GaL(E)$. 
We conclude that $\oline C^{-1}([\1]) = \Aut(A)_E$. 
\end{proof}

\begin{example}
  \label{ex:2.5}
Let $q_\V \: \V \to M$ be a smooth $\K$-vector bundle on the 
compact manifold $M$ and $\Aut(\V)$ the group of smooth bundle isomorphisms. 
Then each element $\phi$ of this group 
permutes the fibers of $\V$, hence induces a diffeomorphism $\phi_M$ of $M$. 
We thus obtain an exact sequence of groups 
$$ \1 \to \Gau(\V) \to \Aut(\V) \to \Diff(M)_{[\V]} \to \1, $$
where $\Gau(\V) = \{ \phi \in \Aut(\V) \: \phi_M = \id_M\}$ 
is the gauge group of $\V$ and 
$$\Diff(M)_{\V} = \{ \psi \in \Diff(M)\: \psi^*\V \cong \V\} $$ 
is the set of all diffeomorphisms $\psi$ of $M$ lifting to 
automorphisms of $\V$ (cf.\ Remark~\ref{rem:2.2}(b)). 
The group $\Diff(M)$ carries a natural Fr\'echet--Lie group structure 
for which  $\Diff(M)_\V$ is an open subgroup, hence also a Lie group. 
Furthermore, it is shown in \cite{ACM89} that $\Aut(\V)$ and $\Gau(\V)$ 
carry natural Lie group structures for which 
$\Aut(\V)$ is a Lie group extension of $\Diff(M)_\V$ by 
$\Gau(\V)$. 

Consider the CIA $A := C^\infty(M,\K)$ and recall from Example~\ref{ex:1.5}(a) that 
the space $E := C^\infty(M,\V)$ of smooth sections of $\V$ is a finitely 
generated projective $A$-module. 
The action of $\Aut(\V)$ on $\V$ induces 
an action on $E$ by 
$$ \phi_E(s)(x) := \phi.s(\phi_M^{-1}(x)). $$
For any smooth function $f \: M \to \K$, we now have 
$\phi_E(fs)(x) := f(\phi_M^{-1}(x))\cdot \phi_E(s)(\phi_M^{-1}(x))$, 
i.e., $\phi_E(fs) = (\phi_M.f)\cdot \phi_E(s)$. We conclude that 
$\Aut(\V)$ acts on $E$ by semilinear automorphisms of $E$ and that we 
obtain a commutative diagram 
$$\begin{matrix}
\Gau(\V) & \to & \Aut(\V) & \to &\Diff(M)_{[\V]} \cr 
\mapdown{}&  & \mapdown{} & &\mapdown{} \cr 
\GL_A(E) & \to & \GaL(E) & \to &\Aut(A)_E.
\end{matrix}$$
Next, we recall that 
\begin{equation}
  \label{eq:2.1}
 \Aut(A) = \Aut(C^\infty(M,\K)) \cong \Diff(M) 
\end{equation}
(cf.\ \cite[Thm.~IX.2.1]{Ne06}, \cite{Bko65}, \cite{Gra06}, \cite{Mr05}). 
Applying \cite[Prop.~4]{Ko76} to the Lie group $G = \Z$, it follows that the 
map $\Aut(\V) \to \GaL(E)$ is a bijective group homomorphism. 
Let us recall the basic idea of the argument. 

First we observe that the vector bundle $\V$ can be reconstructed 
from the $A$-module $E$ as follows. 
For each $m \in M$, we consider the maximal closed ideal 
$I_m := \{ f \in A \: f(m) = 0\}$ and associate the vector space 
$E_m := E/I_m E$. Using the local triviality of the vector bundle $\V$, 
it is easy to see that $E_m \cong \V_m$. We may thus recover 
$\V$ from $E$ as the disjoint union 
$$ \V = \bigcup_{m \in M} E_m. $$

Any $\phi \in \GaL(E)$ defines an automorphism $\phi_A$ of $A$, which 
we identify with a diffeomorphism $\phi_M$ of $M$ via 
$\phi_A(f) := f \circ \phi_M^{-1}$. Then $\phi_A(I_m) = I_{\phi_M(m)}$ 
implies that $\phi$ induces an isomorphism of vector bundles 
$$ \V \to \V, \quad s + I_m E \mapsto \phi(s + I_mE) = \phi(s) + I_{\phi_M(m)}E. $$
Its smoothness follows easily by applying it to a set of 
sections $s_1, \ldots, s_n$ which are linearly independent in $m$. 
This implies that each element $\phi \in \GaL(E)$ corresponds to 
an element of $\Aut(\V)$, so that the vertical arrows in the diagram above 
are in fact isomorphisms of groups. 

Finally, we take a look at the Lie structures on these groups.
A priori, the automorphism group $\Aut(A)$ of a CIA carries no natural 
Lie group structure, but the group isomorphism  $\Diff(M) \to \Aut(A)$ from 
\eqref{eq:2.1}  
defines a smooth action of the Lie group $\Diff(M)$ on $A$. 
Indeed, this can be derived quite directly from the smoothness of the map 
$$ \Diff(M) \times C^\infty(M,\K) \times M \to \K, \quad 
(\phi,f,m) \mapsto f(\phi^{-1}(m)) $$
which is smooth because it is a composition of the smooth action map 
$\Diff(M)\times M \to M$ and the smooth evaluation map 
$A \times M \to M$ (cf.\ \cite[Lemma~A.2]{NeWa08}).

Since the vector bundle $\V$ can be embedded into a trivial 
bundle \break $M \times \K^n$, 
we obtain a topological embedding of $\Gau(\V)$ as a closed subgroup of 
$\Gau(M \times \K^n) \cong C^\infty(M,\GL_n(\K)) \cong \GL_n(A)$. 
Accordingly, we obtain an embedding $E \into A^n$ and 
the preceding discussion yields an identification of 
$\GL_A(E)$ with $\Gau(\V)$ as the same closed subgroups of 
$\GL_n(A)$ (cf.\ Proposition~\ref{prop:1.4}). 
Since both groups are locally exponential Lie groups, 
the homeomorphism $\Gau(\V) \to \GL_A(E)$ is an isomorphism of 
Lie groups (cf.\ \cite[Thm.~IV.1.18]{Ne06} and \cite{GN08} 
for more details 
on locally exponential Lie groups). 
\end{example}

\section{Lie group extensions associated to projective modules} \label{sec:3}

In this section we consider a 
Lie group $G$, acting smoothly by automorphisms on the CIA $A$. 
We write $\mu_A \: G \to \Aut(A)$ for the corresponding 
homomorphism. For each right $A$-module 
$E$,  we then consider the subgroup 
$$ G_E := \{ g \in G \: E^g \cong E \} = \mu_A^{-1}(\Aut(A)_E), $$
where we write $E^g := E^{\mu_A(g)}$ for the corresponding twisted module 
(cf.\ Remark~\ref{rem:2.2}(a)). The main result of this section is 
Theorem~\ref{thm:3.3} which asserts 
that for $G = G_E$, the pull-back of the group extension $\hat\GaL(E)$ of 
$\Aut(A)_E$ by $\GL_A(E)$ yields a Lie group extension $\hat G$ of $G$ by 
$\GL_A(E)$. 

\begin{proposition}
  \label{prop:3.1}
If $E$ is a finitely generated projective right $A$-module, 
then the subgroup $G_E$ of $G$ is open. 
In particular, we have $\mu_A(G) \subeq \Aut(A)_E$ if $G$ is connected. 
\end{proposition}

\begin{proof} 
Since $E$ is finitely generated and projective, it is isomorphic 
to an $A$-module of the form $pA^n$ for some idempotent $p \in M_n(A)$.
We recall from Remark~\ref{rem:2.2}(c) that for any automorphism $\psi \in \Aut(A)$ 
and $\gamma \in \GL_n(A)$ with $M_n(\psi)^{-1}(p) = \gamma^{-1}p\gamma$, 
the maps 
$$ M_n(\psi)^{-1}(p) A^n \to (p A^n)^\psi, x \mapsto \psi^{(n)}(x),\quad 
M_n(\psi)^{-1}(p) A^n \to pA^n, s \mapsto \gamma\cdot s $$
are isomorphisms of $A$-modules. 
According to Proposition~\ref{prop:1.1}, 
all orbits of the group $\GL_n(A)_0$ in 
$\Idem(M_n(A))$ are connected open subsets of \break 
$\Idem(M_n(A))$, hence coincide with its 
connected components.
Therefore the subset 
$\{ g \in G \: g.p \in \GL_n(A)_0.p\}$ of ${G}$ is open. 
In view of Lemma~\ref{lem:1.3}(2), this open subset is contained in 
the subgroup $G_E$, hence $G_E$ is open. 
\end{proof}

From now on we assume that $G = G_E$. 
Then we obtain a group extension 
$$ \1 \to \GL_A(E) \to \hat G  \sssmapright{q} G \to \1, $$
where $q(\phi,g) = g$, and 
\begin{equation}
  \label{eq:3.1}
\hat G := 
\{ (\phi,g) \in \GaL(E) \times G \: (\phi, \mu_A(g)) \in \hat\GaL(E) \} 
\cong \mu_A^*\hat\GaL(E) 
\end{equation}
acts on $E$ via $\pi(\phi,g).s := \phi(s)$ 
by semilinear automorphisms. 
The main result of the present section is that 
$\hat G$ carries a natural Lie group structure and that it is a Lie group 
extension of $G$ by $\GL_A(E)$. Let us make this more precise: 

\begin{definition}
  \label{def:3.2}
An {\it extension of Lie groups} is a short exact sequence 
$$ \1 \to N \sssmapright{\iota} \hat G \sssmapright{q} G \to \1 $$
of Lie group morphisms, for which 
$\hat G$ is a smooth (locally trivial) principal 
$N$-bundle over $G$ with respect to the right action of $N$ given by 
$(\hat g,n) \mapsto \hat gn$. In the following, we identify 
$N$ with the subgroup $\iota(N) \trile \hat G$. 
\end{definition}

\begin{theorem}
  \label{thm:3.3}
If $A$ is a CIA, $G$ is a 
Lie group acting smoothly on $A$ by $\mu_A \: G \to \Aut(A)$, 
and $E$ is a finitely generated projective 
right $A$-module with $\mu_A(G) \subeq \Aut(A)_E$, then 
$\GL_A(E)$ and $\hat G$ carry natural Lie group structures such that the short exact sequence $$ \1 \to \GL_A(E) \to \hat G \sssmapright{q} G \to \1 $$
defines a Lie group extension of $G$ by $\GL_A(E)$. 
\end{theorem}

\begin{proof}
In view of Proposition~\ref{prop:1.4}, 
 $\End_A(E)$ is a CIA and its unit group  
$\GL_A(E)$ is a Lie group. The assumption 
$\mu_A(G) \subeq \Aut(A)_E$ implies that $G = G_E$, so that 
the group $\hat G$ is indeed a group extension of $G$ by $\GL_A(E)$. 

Choose $n$ and $p \in \Idem(M_n(A))$ with $E \cong p A^n$ 
and let $U_p$ be as in Proposition~\ref{prop:1.1}. In the following we identify 
$E$ with $p A^n$.  
We write $g * a := M_n(\mu_A(g))(a)$ for the smooth action of 
$G$ on the CIA $M_n(A)$ and $g \sharp x := \mu_A(g)^{(n)}(x)$ for the 
action of $G$ on $A^n$, induced by the smooth action of $G$ on $A$. 
Then 
$U_G := \{ g \in G \: g * p \in U_p\}$ 
is an open neighborhood of the identity in $G$, and we have a  map 
$$ \gamma \: U_G \to \GL_n(A), \quad g \mapsto s_{g * p} := p\cdot g * p 
+ (\1-p)\cdot (\1- g * p), $$
which, in view of Proposition~\ref{prop:1.1}, satisfies 
\begin{equation}
  \label{eq:3.2}
\gamma(g) (g * p) \gamma(g)^{-1} = p \quad \hbox{ for all } \quad g \in U_G.
\end{equation}
This implies in particular that the natural action of the pair 
$$(\gamma(g),\mu_A(g)) \in \GL_n(A) \rtimes \Aut(A) \cong \GL_A(A^n) \rtimes  
\Aut(A) $$ 
on $A^n$ by $(\gamma(g),\mu_A(g))(x) = \gamma(g)\cdot (g \sharp x)$ 
preserves the submodule $p A^n= E$, and that we thus get a map 
$$ S_E \: U_G \to \GaL(E) = \GaL(pA^n), \quad S_E(g)(s) := \gamma(g)\cdot 
(g \sharp s). $$
For $g \in G$, $s \in E$ and $a \in A$, we then have 
$$ S_E(g)(s.a) = \gamma(g)\cdot (g \sharp (sa)) = \gamma(g)\cdot (g \sharp s) \cdot 
(g * a)
= S_E(g)(s).(g.a), $$
which shows that 
$\sigma \: U_G \to \hat G, g \mapsto (S_E(g),g)$
is a section of the group extension $q \: \hat G \to G$.
We now extend $\sigma$ in an arbitrary fashion 
to a map $\sigma \: G \to \hat G$ with $q\circ \sigma = \id_G$. 
 
Identifying $\GL_A(E)$ with the kernel of the factor map 
$q \: \hat  G \to G$, we obtain for $g,g' \in G$ 
an element 
$$ \omega(g,g') := \sigma(g)\sigma(g') \sigma(gg')^{-1} \in \GL_A(E), $$
and this element is given for $g,g' \in U_G$ by 
$$ \omega(g,g') 
= \gamma(g)\cdot g * \big(\gamma(g') \cdot g^{-1} * \gamma(gg')^{-1}\big) 
= \gamma(g)\cdot (g * \gamma(g')) \cdot \gamma(gg')^{-1}. $$
Since all the maps involved are smooth, the 
preceding formula shows immediately 
that $\omega$ is smooth on the open identity 
neighborhood 
$$\{ (g,g') \in U_G \times U_G 
\: g,g',gg' \in U_G\}$$ of $(\1,\1)$ in $G \times G$. 

Next we observe that the map 
$$ S \: G \to \Aut(\GL_A(E)), 
\quad S(g)(\phi) := \sigma(g)\phi\sigma(g)^{-1}  $$
has the property that the corresponding map 
$$ U_G \times \GL_A(E) \to \GL_A(E), \quad (g,\phi) \mapsto S(g)(\phi) = 
\gamma(g)(g * \phi)\gamma(g)^{-1} $$
is smooth because $\gamma$ and the action of $G$ on $M_n(A)$ are smooth. 

In the terminology of \cite[Def.~I.1]{Ne07}, this means that 
$\omega \in C^2_s(G,\GL_A(E))$ and that $S \in C^1_s(G,\Aut(\GL_A(E)))$. 
We claim that we even have $\omega \in C^2_{ss}(G,\GL_A(E))$, i.e., 
for each $g \in G$, the function 
$$ \omega_g \: G \to \GL_A(E), \quad x \mapsto 
\omega(g,x) \omega(gxg^{-1},g)^{-1}
= \sigma(g)\sigma(x)\sigma(g)^{-1}\sigma(gxg^{-1})^{-1} $$
is smooth in an identity neighborhood of $G$. 

{\bf Case 1:} First 
we consider the case $g \in U_G' := \{h \in G \: h * p \in  \GL_n(A).p\}$. 
The map $\gamma \: U_G \to\GL_n(A)$ extends to a map 
$\gamma \: U_G' \to \GL_n(A)$ satisfying 
$g * p = \gamma(g)^{-1}p \gamma(g)$ for each $g \in U_G'$. Then 
$S_E'(g)(s) := \gamma(g)\cdot (g \sharp s)$ 
defines an element of $\GaL(E)$, and we have 
$\sigma(g) = (\phi(g) S_E'(g),g)$ for some $\phi(g) \in \GL_A(E)$. 
For $x \in U_G \cap g^{-1}U_G g$ we now have 
\begin{align*}
&\ \ \ \omega_g(x)s 
= \sigma(g)\sigma(x)\sigma(g)^{-1}\sigma(gxg^{-1})^{-1}s\cr
&= \phi(g)\gamma(g)\cdot 
g \sharp \Big(\gamma(x) \cdot x\sharp \Big(g^{-1} \sharp \big(\gamma(g)^{-1} 
\phi(g)^{-1} \cdot \big[(gx^{-1}g^{-1}) \sharp \gamma(gxg^{-1})^{-1}s\big]\big)\Big)\Big) \cr
&= \phi(g) \gamma(g)(g  * \gamma(x)) 
\Big((gxg^{-1}) * \big(\gamma(g)^{-1}\cdot \phi(g)^{-1}\big)\Big) 
\cdot \gamma(gxg^{-1})^{-1}s. 
\end{align*}
We conclude that 
$$ \omega_g(x)
= \phi(g) \gamma(g)(g  * \gamma(x)) 
\Big((gxg^{-1}) * \big(\gamma(g)^{-1}\cdot \phi(g)^{-1}\big)\Big) 
\cdot \gamma(gxg^{-1})^{-1}. $$
If $g$ is fixed, all the factors in this product are smooth $\GL_n(A)$-valued 
functions of $x$ in the identity neighborhood $U_G \cap g^{-1}U_Gg$, 
hence $\omega_g$ is smooth on this set. 
 
{\bf Case 2:} Now we consider the case $g * p \not\in \GL_n(A).p$. 
Since $g * p$ corresponds to the right $A$-module $(g * p)A^n \cong 
(pA^n)^{\mu_A(g)^{-1}} 
= E^{\mu_A(g)^{-1}}$ and 
$E^{\mu_A(g)^{-1}} \cong E$ follows from $g \in G_E = G$, 
there exists an element $\eta(g) \in \GL_{2n}(A)$ with 
$\eta(g)(g * \tilde p)\eta(g)^{-1} = \tilde p$ for 
$\tilde p =\pmat{p & 0 \cr 0 & 0 \cr}\in M_{2n}(A)$ 
(Lemma~\ref{lem:1.3}). We now have 
$pA^n \cong \tilde p A^{2n}$, and $g * \tilde p \in \GL_{2n}(A).\tilde p$, 
so that the assumptions of Case $1$ are satisfied with $2n$ instead of $n$. 
Therefore $\omega_g$ is smooth in an identity neighborhood. 

The multiplication in $\hat G = \GL_A(E) \sigma(G)$ is given by the formula 
\begin{equation} \label{eq:3.3}
n\sigma(g) \cdot n'\sigma(g') = \big(n S(g)(n') \omega(g,g')\big) \sigma(gg'),
\end{equation}
so that the preceding arguments imply that 
$(S,\omega)$ is a smooth factor system 
in the sense of \cite[Def.~II.6]{Ne07}. Here the algebraic conditions 
on factor systems follow from the fact that \eqref{eq:3.3} defines a group 
multiplication. We now derive from \cite[Prop.~II.8]{Ne07} 
that $\hat G$ carries a natural Lie group 
structure for which the projection map $q \: \hat G \to G$ defines a 
Lie group extension of $G$ by $\GL_A(E)$. 
\end{proof}

\begin{proposition}
  \label{prop:3.4} {\rm(a)} The group $\hat G$ acts smoothly on $E$ by 
$(\phi,g).s := \phi(s).$

\nin{\rm(b)} Let $p_1 \: \hat G \to \GaL(E)$ denote the projection to the first 
component. For a Lie group $H$, a group homomorphism 
$\Phi \: H \to \hat G$ is smooth if and only if \break 
$q \circ \Phi \: H \to G$ is smooth and the action of $H$ on 
$E$, defined by \break $p_1 \circ \Phi \: H \to \GaL(E)$, is smooth. 

\nin{\rm(c)} The Lie group extension $\hat G$ of $G$ 
splits if and only 
if there is a smooth action of $G$ on $E$ by semilinear automorphisms 
which is compatible with the action of $G$ on $A$ in the sense that 
the corresponding homomorphism $\pi_E \: G \to \GaL(E)$ satisfies 
\begin{equation}
  \label{eq:3.4}
\pi_E(g) \circ \rho_E(a) = \rho_E(\mu_A(g)a) \circ \pi_E(g) \quad
\hbox{ for } \quad g \in G. 
\end{equation}
\end{proposition}

\begin{proof}
(a) Since $\hat G$ acts by semilinear automorphisms of $E$ which 
are continuous and hence smooth, it suffices to see that the action map 
\break $\hat G \times E \to~E$ is smooth on a set of the form $U \times E$, 
where $U \subeq\hat G$ is an open $\1$-neighborhood. 
With the notation of the proof of Theorem~\ref{thm:3.3}, let 
$$ U := \GL_A(E) \cdot \sigma(U_G), $$
which is diffeomorphic to $\GL_A(E) \times U_G$ via the map 
$(\phi,g) \mapsto \phi \cdot\sigma(g)$. Now it remains to observe that 
the map 
$$ \GL_A(E) \times U_G \times E \to E, \quad 
(\phi,g,s) \mapsto \phi(S_E(g)s) = \phi(\gamma(g)\cdot (g\sharp s)) $$
is smooth, which follows from the smoothness of the action of 
$\GL_A(E)$ on $E$, the smoothness of $\gamma \: U_G \to \GL_n(A)$ and the 
smoothness of the action of $G$ on $A^n$. 

(b) If $\Phi$ is smooth, then 
$q \circ \Phi$ is smooth and (a) implies 
that $f := p_1 \circ \Phi$ defines a smooth action of $H$ on $E$. Suppose, 
conversely, that $q \circ \Phi$ is smooth and that $f$ 
defines a smooth action on $E$. Let $U_G$ and $U$ be as in (a) and 
put $W := \Phi^{-1}(U)$. Since $q \circ \Phi$ is continuous, 
$$ W = \Phi^{-1}(q^{-1}(U_G)) = (q \circ \Phi)^{-1}(U_G) $$
is an open subset of $H$. Since $\Phi$ is a group homomorphism, it suffices
to verify its smoothness on $W$. We know from (a) that 
the map 
$$ \tilde S_E \: U_G \times E \to E, \quad (g,s) \mapsto S_E(g)s $$
is smooth. For $h \in W$ we have 
$\Phi(h) = f_1(h) \sigma(f_2(h))$, where 
$f_1(h) \in \GL_A(E)$ and 
$f_2 := q \circ \Phi \: W \to U_G$ is a smooth map. Therefore the map 
$$ W \times E \to E, \quad (h,s) \mapsto f(h)\big(\sigma(f_2(h))^{-1}.s\big) 
= f(h) S_E(f_2(h))^{-1}.s = f_1(h).s $$
is smooth. If $e_1, \ldots, e_n \in A^n$ denote the canonical basis elements 
of the right $A$-module $A^n$, then we conclude that all maps 
$$ W \to \GL_A(E), \quad h \mapsto f_1(h) \cdot e_i 
= f_1(h)p e_i $$ 
are smooth because $p e_i \in E$. Hence all columns of the matrix 
$f_1(h)$ depend smoothly on $h$, and thus $f_1 \: W \to \GL_A(E) \subeq M_n(A)$ 
is smooth.  This in turn implies that 
$\Phi(h) = f_1(h) \sigma(f_2(h))$ is smooth on $W$, hence on $H$ because 
it is a group homomorphism. 

(c) First we note that any homomorphism 
$\Phi \: G \to \hat G$ is of the form $\tilde\sigma(g) = (f(g),g),$
where $f \: G \to \GaL(E)$ is a homomorphism satisfying 
$(f(g),\mu_A(g)) \in \hat\GaL(E)$. 

If the extension $\hat G$ of $G$ by $\GL_A(E)$ 
splits, then there is such a smooth $\Phi$, and then (a) implies that 
$\pi_E = f$ defines a smooth action of $G$ on $E$, satisfying all 
requirements. 

If, conversely, $\pi_E \: G \to \GaL(E)$ defines a smooth action with 
\eqref{eq:3.4}, then 
the map $\Phi = (\pi_E, \id_G) \: G \to \hat G$ is a group homomorphism whose 
smoothness follows from (b), and therefore the Lie group 
extension $\hat G$ splits.
\end{proof}

\begin{examples}
  \label{ex:3.5}
(a) If $A$ is a Banach algebra, then $\Aut(A)$ carries a natural 
Banach--Lie group structure (cf.\ \cite{HK77}, \cite[Prop.~IV.14]{Ne04}). 
For each finitely generated projective module $pA^n$, $p \in \Idem(M_n(A))$, 
the subgroup $G := \Aut(A)_E$ is open (Proposition~\ref{prop:3.1}), and we thus obtain a Lie group 
extension 
$$ \1\to \GL_A(E) \into \hat G \onto G = \Aut(A)_E \to \1.$$

(b) For each CIA A, the Lie group $G := A^\times$ acts smoothly by conjugation on $A$ 
and $g \mapsto \rho_E(g^{-1})$ defines a smooth action of $G$ on $E$ 
by semilinear automorphisms. This leads to a homomorphism 
$$ \sigma \: A^\times \to \hat G = \{ (\phi,g) \in \GaL(E) \times G 
\: (\phi, \mu_A(g)) \in \hat\GaL(E) \},\ 
g \mapsto (\rho_E(g^{-1}),g), $$
splitting the Lie group extension $\hat G$ (Proposition~\ref{prop:3.4}). 

Note that for any CIA $A$, we have $Z(A)^\times = Z(A^\times)$ because 
$A^\times$ is an open subset of $A$, so that its centralizer coincides with 
the center $Z(A)$ of $A$. We also note that 
$\rho_E(Z(A)) \subeq \End_A(E)$ and that 
the direct product group $\GL_A(E) \times A^\times$ 
acts on $E$ by $(\phi,g).s := \phi \circ \rho_E(g^{-1})s$, 
where the pairs $(\rho_E(z), z^{-1})$, $z \in Z(A^\times)$, act 
trivially. 

If, in addition, $A$ is Mackey complete, 
the Lie group $\GL_A(E) \times A^\times$ and both factors are locally 
exponential, the subgroup 
$$\Delta_Z := \{ (\rho_E(z),z^{-1})\: z \in Z(A^\times)\}$$ 
is a central Lie subgroup and the Quotient Theorem in \cite{GN08}  
(see also \cite[Thm.~IV.2.9]{Ne06}) 
implies that $(\GL_A(E) \times A^\times)/Z(A^\times)$ carries a locally 
exponential Lie group structure. 

If, in addition, $E$ is a faithful $A$-module, then $\Delta_Z$ coincides with 
the kernel of the action of $\GL_A(E) \times A^\times$ on $E$, so that 
the Lie group \break $(\GL_A(E) \times A^\times)/Z(A^\times)$ injects into $\GaL(E)$. 
If, moreover, all automorphisms of $A$ are inner, we have 
$\Aut(A) \cong A^\times/Z(A^\times),$
which carries a locally exponential Lie group structure 
(\cite[Thm.~IV.3.8]{Ne06}). We obtain 
$$ \GaL(E) \cong (\GL_A(E) \times A^\times)/Z(A^\times), $$
and a Lie group extension 
$$ \1 \to \GL_A(E) \to \GaL(E) \onto \Aut(A) \to\1. $$

(c) (Free modules) Let $n \in \N$ and let $p = \1 \in M_n(A)$. 
Then $pA^n = A^n$, $\GL_A(E) \cong \GL_n(A)$ acts by left multiplication, and 
$$ \GaL(E) \cong \GL_n(A) \rtimes \Aut(A) $$
is a split extension. For any smooth Lie group action 
$\mu_A \: G \to \Aut(A)$, we accordingly get $G_E = G$ and a split extension 
$\hat G \cong \GL_n(A) \rtimes G.$

(d) Let $A = B(X)$ denote the Banach algebra of all bounded operators 
on the complex Banach space $X$. 
If $p \in A$ is a rank-$1$-projection, 
$$ pA \cong X' = \Hom(X,\C) $$
is the dual space, considered as a right $A$-module, the module structure 
given by 
$\phi.a := \phi \circ a.$
In this case $pAp \cong \C$, $\GL_A(E) \cong\C^\times$, and for the group 
$G := \PGL(X) := \GL(X)/\C^\times$, 
acting by conjugation on $A$, we obtain the central extension 
$$ \1 \to \C^\times \to \hat G \cong \GL(X) \to G = \PGL(X) \to \1. $$
\end{examples}

\section{The corresponding Lie algebra extension} \label{sec:4}

We now determine the Lie algebra of the Lie group 
$\hat G$ constructed in Theorem~\ref{thm:3.3}. This will lead us from semilinear automorphisms 
of a module to derivative endomorphisms. The relations to 
connections in the context of non-commutative geometry will be discussed in 
Section~\ref{sec:5} below.  

\begin{definition}
  \label{def:4.1}
We write $\gl_A(E)$ 
for the Lie algebra underlying the associative 
algebra $\End_A(E)$ and 
\begin{align*}
&\dend(E) \\ 
&:= \{ \phi \in \End_\K(E) \: (\exists D_\phi \in \der(A))(\forall a \in A)\ \ 
[\phi,  \rho_E(a)] = \rho_E(D_\phi(a))\}. 
\end{align*}
for the Lie algebra of {\it derivative endomorphisms of $E$} 
(cf.\ \cite{Ko76}). 
We write 
$$ \hat\dend(E) := \{ (\phi,D) \in \End_\K(E)\times \der(A) \: 
(\forall a \in A)\ [\phi, \rho_E(a)] = \rho_E(D.a)\}. $$
We then have a short exact sequence 
$$ \0 \to \gl_A(E) \to \hat\dend(E) \to \der(A)_E \to \0 $$
of Lie algebras, where 
$$\der(A)_E = \{ D \in \der(A) \: (\exists \phi \in \dend(E))\ D_\phi = D\} $$
is the image of the homomorphism $\dend(E) \to \der(A)$. 
\end{definition}

\begin{example}
  \label{ex:4.2}
If $E = C^\infty(M,\V)$ is the space of smooth sections of 
the vector bundle $\V$ with typical fiber $V$ 
on the compact manifold $M$, then it is interesting to 
identify the Lie algebra $\dend(E)$. 
From the short exact sequence 
$$ \0 \to \gl_A(E) = C^\infty(\End(\V)) \into \dend(E) 
\onto {\cal V}(M) = \der(A) \to \0, $$
it easily follows that the Lie algebra $\dend(E)$ can be 
identified with the Lie algebra ${\cal V}(\Fr \V)^{\GL(V)}$ 
of $\GL(V)$-invariant vector fields on the frame bundle $\Fr \V$ 
(cf.\ \cite{Ko76} for details).  
\end{example}

\begin{lemma}
  \label{lem:4.3}
For each $a \in A$ we have 
$(\rho_E(a), -\ad a) \in \hat\dend(E)$ and in particular 
$\rho_E(A) \subeq \dend(E)$. 
\end{lemma}

\begin{proof}
For each $b \in A$ we have 
$\rho_E(a)(s.b) - \rho_E(a)(s).b = s.(ba-ab)= s.(-\ad a(b)).$ 
\end{proof}

\begin{lemma}
  \label{lem:4.4}
Let $p \in \Idem(M_n(A))$ and $E := pA^n$. We define 
$$\gamma \: \der(A) \to M_n(A), \quad \gamma(D) := (2p-\1)\cdot (D.p). $$
Then $[p, \gamma(D)] = D.p$ and the operator 
$$ \nabla_Ds := \gamma(D)s + D.s $$ 
on $A^n$ preserves $E = pA^n$ and $(\nabla_D,D) \in \hat\dend(E)$. 
\end{lemma}

\begin{proof}
From $p^2 = p$ we immediately get 
$D.p = D.p^2 = p\cdot (D.p) + (D.p) \cdot p$, showing that 
$p\cdot (D.p) = p \cdot (D.p) + p \cdot (D.p) \cdot p$, and therefore 
$p \cdot (D.p) \cdot p = 0.$
We  likewise obtain 
$(\1-p) \cdot (D.p) \cdot (\1-p) = 0,$
so that 
$D.p \in pA(\1-p) + (\1-p) A p.$
This leads to 
\begin{align*}
[\gamma(D),p] 
&= (2p-\1)\cdot (D.p) \cdot p - p (2p-\1)\cdot D.p \\
&= - (D.p) \cdot p - p\cdot D.p = - D.(p^2) = - D.p.
\end{align*}

For any $s \in pA^n$ we have $ps = s$ and therefore  
\begin{align*}
p(\gamma(D)s+D.s) 
&= [p,\gamma(D)]s + \gamma(D)ps + p(D.s) 
= (D.p)s + \gamma(D)s + p(D.s) \cr
&= D.(ps) + \gamma(D)s 
= D.s + \gamma(D)s. 
\end{align*}
This implies that $\nabla_D.s \in p A^n = \{ x \in A^n \: px = x\}$.

The remaining assertion follows from the fact that left and right multiplications 
commute and $D.(s.a) = (D.s).a + s.(D.a)$. 
\end{proof}

\begin{remark}
  \label{rem:4.5}
To calculate the Lie algebra $\L(G)$ of a Lie group $G$, 
one may use a local chart $\phi \: U \to \L(G)$, $U \subeq G$ an open $\1$-neighborhood 
and $\phi(\1)=~0$, and consider the Taylor 
expansion of order $2$ of the multiplication 
$$ x * y := \phi(\phi^{-1}(x)\phi^{-1}(y)) = x + y + b_\g(x,y) + \cdots. $$
Then the Lie bracket in $\L(G)$ satisfies 
$$[x,y] = b_\g(x,y) - b_\g(y,x) $$
(cf.\ \cite{Mil84}, \cite{GN08}). 

Suppose that a Lie group extension 
$N \into \hat G \onto G$ is given by a pair $(S,\omega)$ via the multiplication 
on the product set $N \times G$: 
$$ (n,g)(n',g') = \big(n S(g)(n') \omega(g,g'), gg'), $$
where the maps 
$S \: G \to \Aut(N)$ and $\omega \: G \times G \to N$ satisfy: 
\begin{description}
\item[(a)] the map $G \times N \to N, (g,n) \mapsto S(g)n$ is smooth on a set 
of the form $U \times N$, where $U$ is an identity neighborhood of $G$. 
\item[(b)] $\omega$ is smooth in an identity neighborhood with 
$\omega(g,\1) = \omega(\1,g) = \1$. 
\end{description}

Let $\g = \L(G)$,  $\n= \L(N)$ and $\hat\g = \L(\hat G)$ be the corresponding 
Lie algebras. Then we may use a product chart of $\hat G$ in some sufficiently small 
$\1$-neighborhood. Write $\L(S(g)) \in \Aut(\n)$ for the 
Lie algebra automorphism induced by $S(g) \in \Aut(N)$ and put  
$\L(S)(g,v) := \L(S(g))(v)$. 
Then we define $DS \: \g \to \gl(\n)$ by 
$$ DS(x)(v) := T_{(\1,v)}(\L(S))(x,0). $$
We further define 
$$D\omega(y,y') := d^2\omega(\1,\1)((y,0),(0,y')) - d^2\omega(\1,\1)((y',0),(0,y)) $$
and note that this is well-defined, i.e., independent of the chart, because 
(b) implies that $\omega$ vanishes of order $1$ in $(\1,\1)$. 
We thus obtain the second order Taylor expansion 
\begin{align*}
(x,y)*(x',y') 
&= (x + DS(y)(x') + b_\n(x,x') + d^2\omega(\1,\1)((y,0),(0,y')) + \cdots, \\
&\qquad\qquad\qquad\qquad y + y' + b_\g(y,y') + \cdots), 
\end{align*}
which provides the following formula for the Lie bracket in 
$\hat\g$, written as  the product set $\n \times \g$: 
$$ [(x,y),(x',y')] = ([x,x'] + DS(y)(x') - DS(y')(x) + D\omega(y,y'), [y,y']).
$$
\end{remark}

To show that the Lie algebra of the group $\hat G$,  constructed in 
Theorem~\ref{thm:3.3}, 
is the corresponding pull-back on the Lie algebra level, we 
need the following lemma: 

\begin{lemma}
  \label{lem:4.6}
Let $G$ and $H$ be Lie groups with Lie algebra $\g$, resp., $\h$, 
$U \subeq G$ an open $\1$-neighborhood and $\sigma \: U \to H$ a smooth 
map with $\sigma(\1) = \1$. For the map 
$$ \omega \: U \times U \to H, \quad (g,g') \mapsto \sigma(g)\sigma(g')\sigma(gg')^{-1}, $$
we then have 
\begin{align*}
& d^2\omega(\1,\1)((x,0),(0,x')) - d^2\omega(\1,\1)((x',0),(0,x)) \\
& = [T_\1(\sigma)x, T_\1(\sigma)x'] - T_\1(\sigma)[x,x']. 
\end{align*}
\end{lemma}

 \begin{proof}
For the function $\omega_g(g') := \omega(g,g')$, we directly obtain with 
respect to the group structure on the tangent bundle $TH$:  
\begin{align*}
T_\1(\omega_g)x' 
&= \sigma(g)\cdot T_\1(\sigma)x' \cdot \sigma(g)^{-1} 
+ \sigma(g)\big(-\sigma(g)^{-1}T_g(\sigma)(g.x') \sigma(g)^{-1}\big)\cr
&= \sigma(g)\cdot T_\1(\sigma)x' \cdot \sigma(g)^{-1} 
- T_g(\sigma)(g.x') \sigma(g)^{-1}\cr
&= \Ad(\sigma(g)) T_\1(\sigma)x' - \delta^r(\sigma)(x'_l)(g), 
\end{align*}
where $\delta^r(\sigma)\in \Omega^1(G,\h)$ is the right logarithmic derivative 
of $\sigma$ and $x'_l(g) = g.x'$ is the left invariant vector field on $G$, 
corresponding 
to $x'$. Taking derivatives in $g = \1$ with respect to $x$, 
this in turn leads to 
$$ (d\omega)(\1,\1)(x,x') = [T_\1(\sigma)x,T_\1(\sigma)x'] 
- x_l\big(\delta^r(\sigma)(x'_l)\big)(\1). $$
Using the Maurer--Cartan equation (\cite{KM97}) 
$$\dd\delta^r(\sigma)(X,Y) = [\delta^r(\sigma)(X),\delta^r(\sigma)(Y)],$$
we now obtain 
\begin{align*}
&\ \ \ \ d^2\omega(\1,\1)((x,0),(0,x')) - d^2\omega(\1,\1)((x',0),(0,x)) \cr
&= 2[T_\1(\sigma)x,T_\1(\sigma)x'] 
- x_l\big(\delta^r(\sigma)(x'_l)\big)(\1)  
+ x'_l\big(\delta^r(\sigma)(x_l)\big)(\1)  \cr
&= 2[T_\1(\sigma)x,T_\1(\sigma)x'] - \dd\delta^r(\sigma)(x_l, x'_l)(\1) 
- \delta^r(\sigma)([x_l, x'_l])(\1)  \cr
&= 2[T_\1(\sigma)x,T_\1(\sigma)x'] - [\delta^r(\sigma)(x_l), \delta^r(\sigma)(x'_l)](\1) 
- T_\1(\sigma)([x,x']) \cr
&= 2[T_\1(\sigma)x,T_\1(\sigma)x'] 
- [T_\1(\sigma)x,T_\1(\sigma)x']  - T_\1(\sigma)([x,x']) \cr
&= [T_\1(\sigma)x,T_\1(\sigma)x'] - T_\1(\sigma)([x,x']). 
\end{align*}
\end{proof}

Since, in general,  
$\Aut(A)$ does not carry a natural Lie group structure, the Lie algebra 
$\der(A)$ is not literally the Lie algebra of $\Aut(A)$, 
but $\mu_A$ leads to a homomorphism $\L(\mu_A) \: \g \to \der(A)$ 
of Lie algebras, given by 
$$ x.a := \L(\mu_A)(x)(a) := (T\mu_A)(\1,a)(x,0). $$
(cf.\ \cite[App.~E]{GN08} for a 
discussion of these subtle points in the infinite-dimensional context; 
see also \cite[Remark~II.3.6(a)]{Ne06}).

\begin{proposition}
  \label{prop:4.7}
The Lie algebra $\hat\g := \L(\hat G)$ of the Lie group 
$\hat G$ from {\rm Theorem~\ref{thm:3.3}} is isomorphic to 
$$ \L(\mu_A)^*\hat\dend(E)\cong 
\{ (x,\phi) \in \g \times \dend(E) \: (\phi, \L(\mu_A)x) \in 
\hat\dend(E)\}. $$
\end{proposition}

\begin{proof}
First we note that 
\begin{align*}
&\L(\mu_A)^*\hat\dend(E) 
= \{ (x, (\phi, D)) \in \g \times \hat\dend(E) \: \L(\mu_A)x = D \} \cr
&= \{ (x, (\phi, \L(\mu_A)x)) \in \g \times \dend(E) \times \der(A) \: 
(\phi,\L(\mu_A)x) \in \hat\dend(E) \} \cr
&\cong \tilde\g := \{ (x, \phi)\in \g \times \dend(E) \: (\phi,\L(\mu_A)x) \in \hat\dend(E) \}.
\end{align*}
 
Recall from the proof of Theorem~\ref{thm:3.3} the maps 
$$ \gamma(g) = p\cdot (g.p) + (\1-p)\cdot (\1-(g.p)) 
\quad \hbox{ and } \quad 
S_E(g)(s) := \gamma(g)\cdot (g.s). $$
Taking derivatives, we get 
$$ \dot\gamma(x) := T_\1(\gamma)(x) = (2p-\1) \L(\mu_A)(x).p 
= (2p-\1) \cdot (x.p) $$ 
and with Lemma~\ref{lem:4.4} we obtain 
$[\dot\gamma(x),p] = - x.p$. 
We further derive 
$$ T_\1(S_E)(x).s := \dot\gamma(x)\cdot s + x.s \in E, $$ 
and the linear map $T_\1(S_E) \: \g \to \dend(E)$ satisfies 
$$ [T_\1(S_E)(x),\rho_E(a)] = \rho_E(x.a) $$
for each $x \in \g$, which means that 
$(T_1(S_E)(x), x) \in \tilde\g$. 
Now 
$$ \Gamma \: \hat\g = \gl_A(E) \oplus \g \to \tilde\g, \quad 
(\phi,x) \mapsto (\phi + T_\1(S_E)(x), x) $$
is a linear isomorphism with 
\begin{align*}
& [\Gamma(\phi,x), \Gamma(\phi',x')] \\
&= ([\phi,\phi'] + [T_\1(S_E)(x),\phi'] - [T_\1(S_E)(x'),\phi] 
+ [T_\1(S_E)(x),T_\1(S_E)(x')], [x,x']). 
\end{align*}
On the other hand, we have seen in Remark~\ref{rem:4.5} that the Lie 
bracket in $\hat\g = \gl_A(E) \times \g$ is given by 
$$ [(\phi,x),(\phi,x')] = ([\phi,\phi'] + DS(x)(\phi') - DS(x')(\phi) 
+ D\omega(x,x'), [x,x']). $$
From $S(g)(\phi) = \Ad(\gamma(g))(g.\phi)$ we derive that 
$$ DS(x)(\phi) = [\dot\gamma(x),\phi] + x.\phi = [T_\1(S_E)(x), \phi]. $$
To show that $\Gamma$ is an isomorphism of Lie algebras, 
it therefore suffices to show that 
\begin{equation}
  \label{eq:4.1}
[T_\1(S_E)(x),T_\1(S_E)(x')] = D\omega(x,x') + T_\1(S_E)([x,x']) 
\end{equation}
holds for $x,x' \in \g$. 

To verify this relation, we first observe that the smoothness of the 
action of $G$ on $\GL_n(A)$ implies that $\GL_n(A) \rtimes G$ is a Lie group, 
acting smoothly on $A^n$ by $(a,g).s := a(g \sharp s)$. 
Then we consider the smooth map 
$$ \tilde S_E \: U_G \to \GL_n(A) \rtimes G, \quad g \mapsto (\gamma(g),g), $$
also satisfying 
$\omega(g,g') = \tilde S_E(g)\tilde S_E(g')\tilde S_E(gg')^{-1}$
and 
\begin{equation}
  \label{eq:4.2}
\tilde S_E(g).s = S_E(g).s \quad \hbox{ for } \quad s \in E = pA^n.
\end{equation}
Lemma~\ref{lem:4.6} provides the identity 
$$ D\omega(x,x') 
= [T_\1(\tilde S_E)x, T_\1(\tilde S_E)x'] - T_\1(\tilde S_E)[x,x'], $$
in $\gl_A(E)$, so that \eqref{eq:4.2} leads to 
$$ D\omega(x,x') 
= [T_\1(S_E)x, T_\1(S_E)x'] - T_\1(S_E)[x,x'], $$
as linear operators on $E$, and this is \eqref{eq:4.1}. 
\end{proof}

The following general lemma prepares the discussion in 
Example~\ref{ex:4.9} 
below, which exhibits the Lie group $\Aut(\V)$ of automorphisms 
of a vector bundle as one of the Lie group extensions from 
Theorem~\ref{thm:3.3}. 

\begin{lemma}
  \label{lem:4.8}
If $\phi \: G \to H$ is a bijective 
morphism of regular connected Lie groups and 
$\L(\phi) \: \L(G) \to \L(H)$ is an isomorphism of locally convex 
Lie algebras, then $\phi$ is an isomorphism of Lie groups. 
\end{lemma}

\begin{proof}
Let $q_G \: \tilde G \to G$ and $q_H \: \tilde H \to H$ denote 
simply connected universal covering groups with
$\L(q_G) = \id_{\L(G)}$ and $\L(q_H) = \id_{\L(H)}$. 
Then the induced morphims $\tilde\phi \: \tilde G \to \tilde H$ 
is the unique morphism of Lie groups with 
$\L(\tilde \phi) = \L(\phi)$, hence an isomorphism, whose inverse is the 
unique morphism $\psi \: \tilde H \to \tilde G$ with 
$\L(\psi) = \L(\phi)^{-1}$ (\cite[Thm.~IV.1.19]{Ne06}). 
Since $\tilde\phi$ is an isomorphism, $\phi$ is a local 
isomorphism of Lie groups, so that its bijectivity implies that it is 
an isomorphism. 
\end{proof}

\begin{example}
  \label{ex:4.9}
We continue the discussion of Example~\ref{ex:2.5} in the 
light of Theorem~\ref{thm:3.3}. 
Recall that $q_\V \: \V \to M$ denotes a smooth $\K$-vector bundle on the 
compact manifold $M$ and $\Aut(\V)$ its group of smooth bundle isomorphisms. 

Let $G := \Diff(M)_{[\V]}$ and recall that this group acts smoothly 
on the CIA $A = C^\infty(M,\K)$, preserving the equivalence class of the 
projective module $E = \Gamma\V$. 
Let $\hat G$ be the Lie group extension of $G$ by $\GL_A(E)$ from 
Theorem~\ref{thm:3.3}. In view of Proposition~\ref{prop:2.4}, we 
have a smooth representation $\pi \: \hat G \to \GaL(E)$ of 
$\hat G$ on $E$ whose range is a subgroup of 
$\GaL(E)$ containing $\GL_A(E)$ and projecting onto $\Diff(M)_{[\V]} 
\cong \Aut(A)_E$, which implies that the representation 
$\pi$ is a bijection. 
In Example~\ref{ex:2.5} we have seen that $\hat G \cong \Aut(\V)$ as abstract 
groups. 

Next we recall that $G$ is a regular Lie group (\cite[Thm.~38.6]{KM97}), 
and that we have 
seen in Example~\ref{ex:2.5} that $\GL_A(E) \cong \Gau(\V)$ as Lie groups, 
which implies that $\GL_A(E)$ is regular (\cite[Thm.~38.6]{KM97}). Hence 
$\hat G$ is an extension of a regular Lie group by a regular Lie group 
and therefore regular (cf.\ \cite{KM97}, \cite{GN08}). For similar reasons, 
the Lie group $\Aut(\V)$ is regular. 
To see that $\Aut(\V)\cong \hat G$ as Lie groups, it therefore 
suffices to show that the canonical isomorphism 
$\phi \: \Aut(\V) \to \GaL(E) \cong \hat G$ is smooth and that 
$\L(\phi)$ is an isomorphism of topological Lie algebras (Lemma~\ref{lem:4.8}). 

In view of Proposition~\ref{prop:3.4}(b), the smoothness of $\phi$ follows from the 
smoothness of the action of $\Aut(\V)$ on $\Gamma\V$. 
To verify this smoothness, let $P = \Fr \V$ denote the frame bundle
 of $\V$ and recall that 
$$ \Gamma \V  \cong \{ f \in C^\infty(P,V) \: (\forall p \in P)
(\forall k \in \GL(V))\ f(p.k) = k^{-1} f(p)\}. $$
Since the evaluation map of $C^\infty(P,V)$ is smooth, the smoothness 
of the action of $\Aut(\V) \cong \Aut(P)$ on $\Gamma\V$ now follows from 
the smoothness of the action of $\Aut(P)$ on $P$ 
(cf.\ \cite{ACM89}). 

To see that $\hat\g \cong \L(\Aut(\V)) \cong {\cal V}(\Fr \V)^{\GL(V)}$, 
we first recall from Example~\ref{ex:4.2} that $\dend(\V) \cong 
{\cal V}(\Fr \V)^{\GL(V)}$. Hence both $\hat\g$ and $\L(\Aut(\V))$ are 
Fr\'echet--Lie algebras and $\L(\phi)$ is a continuous homomorphism 
inducing bijections $\Gamma(\End(\V)) \to \gl_A(E)$ and 
${\cal V}(M) \to \g$. Hence $\L(\phi)$ is bijective, and the Open Mapping 
Theorem (\cite[Thm.~2.11]{Ru73}) 
implies that it is an isomorphism of topological 
Lie algebras. This completes the proof that $\hat G \cong \Aut(\V)$ 
as Lie groups. 
\end{example}

\begin{remark}
  \label{rem:4.10}
That the extension $\hat\g$ of $\g$ by $\gl_A(E)$ 
splits is equivalent to the existence of a continuous linear map 
$\alpha \: \g \to \gl_A(E)$ for which 
$$ T_\1(S_E) + \alpha  \: \g \to \dend(E) $$ 
is a homomorphism of Lie algebras (cf.\ Proposition~\ref{prop:3.4} 
for the corresponding 
group analog). This is equivalent to the 
existence of a $\g$-module structure on $E$, lifting the action of 
$\g$ on $A$, given by $\L(\mu_A)$. 

If such a homomorphism exists and the group $\hat G$ is regular, then 
there exists a morphism of Lie groups 
$\tilde G_0 \to \hat G$, splitting the pull-back extension 
$q_G^*\hat G$ of $\tilde G_0$ by $\GL_A(E)$. 
\end{remark}

\section{Covariant derivatives} \label{sec:5}

In this short final section, we briefly explain the connections 
between linear splittings of the Lie algebra extension from 
Proposition~\ref{prop:4.7} and covariant derivatives, resp., 
connections, as they occur in non-commutative geometry. 

\begin{definition}
  \label{def:5.1}

(a) There are many ways to construct ``differential forms'' for a 
non-commutative algebra (cf.\ \cite{Co94}). 
One approach, which is closest to our construction, 
is the one described by Dubois-Violette in \cite{DV91} (cf.\ 
\cite{DV88} and \cite{DVM94}): 
First, one considers $A$ as a module of  the Lie algebra $\der(A)$, 
and since the multiplication on $A$ is $\der(A)$-invariant, 
the algebra multiplication provides on the 
Chevalley--Eilenberg complex 
$(C(\der A,A),\dd_A)$ the structure of a differential graded 
algebra. The differential subalgebra 
generated by $A \cong C^0(\der A,A)$ 
and $\dd_A(A)$ is denoted $\Omega_D(A)$. From the inclusion 
of $A$ as a subalgebra, $\Omega_D(A)$ inherits a natural 
$A$-bimodule structure. In particular, $E \otimes_A \Omega^1_D(A)$ is 
defined and a right $A$-module. 

(b) Let $E$ be a right $A$-module. 
A {\it connection} on $E$ is a linear map 
$\nabla \: E \to E \otimes_A \Omega^1_D(A)$
satisfying 
\begin{equation}
  \label{eq:4.1b}
\nabla(s.a) = \nabla(s)a + s \otimes \dd_A a 
\quad \hbox{ for } \quad s \in E, a \in A. 
\end{equation}
Since the elements of $\Omega^1_D(A)$ are linear maps $\der(A) \to A$, 
each element of $E \otimes_A \Omega^1_D(A)$ defines a
linear map $\der(A) \to E$. If $i_D \: \Omega^1_D(A) \to A, 
\alpha \mapsto \alpha(D),$ 
denotes the evaluation map, we thus obtain 
for each derivation $D \in \der A$ 
a linear map, the corresponding {\it covariant derivative}, 
$$ \nabla_D := (\id_E \otimes i_D) \circ \nabla\: E \to E, $$
where we identify $E \otimes_A A$ with $E$. The covariant derivative satisfies 
$$ \nabla_D(sa) = \nabla_D(s)a + sDa \quad \hbox{ for } \quad s \in E, 
a \in A. $$
\end{definition}

\begin{remark}
  \label{rem:5.2}
(a) For any connection $\nabla$ 
and $D \in \der(A)$, we have $(\nabla_D, D) \in \hat\dend(E)$. 
In particular, we have $\der(A)_E = \der A$ whenever a connection exists,  
and in this case any connection $\nabla$ defines a splitting 
of the Lie algebra extension 
$\hat\dend(E) \to \der(A)$ of $\der(A)$ by $\gl_A(E)$ 
(cf.\ Definition~\ref{def:4.1}). 
In this sense, we call any linear section of this Lie algebra extension 
a {\it covariant derivative on $E$}. 

(b) (Covariant coordinates) We have already seen in Lemma~\ref{lem:4.3} 
that 
for each $a \in A$, the operator $\rho_E(a)$ is contained in $\dend(E)$ 
and satisfies $D_{\rho_E(a)} = - \ad a$. If $\nabla$ is a connection, we therefore have 
$$ \hat\rho_E(a) := \rho_E(a) + \nabla_{\ad a} \in \gl_A(E),
\quad \mbox{ i.e., } \quad  [\hat \rho_E(A),\rho_E(A)] = \{0\}. $$
In the context of non-commutative geometry, the 
operators $\hat\rho_E(a)$ are called {\it covariant coordinates} 
because they commute with all ``coordinate operators'' 
$\rho_E(a)$, $a \in A$ (cf.\ \cite{Sch01}, \cite{JSW01}). 
\end{remark}

\begin{remark}
  \label{rem:5.3}
(a) If a $E$ 
is a finitely generated projective module, then it is of the form 
$p A^n$ for some idempotent $p \in M_n(A)$. In this case we have 
the {\it Levi--Civita connection}, given by 
$$ \nabla(s):= p \cdot \dd_{A^n}(s), $$
where 
$\dd_{A^n} \: A^n  \to A^n \otimes_A \Omega^1_D(A) \cong \Omega^1_D(A)^n, 
(a_i) \mapsto (\dd_A(a_i))$ 
is the canonical connection of the free right module~$A^n$. 
Note that the Levi--Civita connection 
is not an intrinsic object, it depends on the 
embedding $E \into A^n$ and the module 
complement, all of which is encoded in the choice of the idempotent~$p$. 
For a derivation $D \in \der(A)$, the operator $\nabla_D \in \dend(A^n)$, 
defined  by 
$$ \nabla_D(s) = p \cdot D.s $$
is the covariant derivative corresponding to the 
Levi--Civita connection. 

In Lemma~\ref{lem:4.4} we have seen that 
$\nabla_D' s := (2p-\1)(D.p)s + D.s$ 
also defines a covariant derivative on $E$. 
In view of $ps = s$ for $s \in E$, we have 
$p(D.p)s = p(D.p)(ps) = 0$ (see the proof of Lemma~\ref{lem:4.4}), so that 
$$ \nabla_D' s = -(D.p)s + D.s = D.(ps) - (D.p)s = p (D.s) = \nabla_D s. $$ 

(b) For $E = p A^n$ as above, any connection $\nabla'$ on $E$ is of the form 
$$ \nabla_\alpha = p \dd_{A^n} + \alpha, $$
where $\alpha \in M_n(\Omega^1_D(A))$ satisfies $\alpha = p \alpha p$, i.e., 
$\alpha \in p M_n(\Omega^1_D(A))p$. 
Here we use that for any other connection $\nabla'$ we have 
\begin{align*}
 \nabla'- \nabla &\in 
\Hom(E, E \otimes_A \Omega^1_D(A)) 
= \Hom(pA^n, pA^n \otimes_A \Omega^1_D(A)) \\
&\cong p\Hom(A^n, A^n \otimes_A \Omega^1_D(A))p 
= p M_n(\Omega^1_D(A))p.
\end{align*}

For any gauge transformation $g \in \GL_A(E)$, we then  have 
$$ \nabla(g.s) = p  \dd_{A^n}(g.s) 
= p( \dd_{M_n(A)} g \cdot s + g.\dd_{A^n}s) 
= g.\big(\nabla(s) + g^{-1}\cdot \dd_{M_n(A)}(g) \cdot s\big), $$
which for $\nabla^g(s) := g^{-1}\nabla(g.s)$ and the left logarithmic 
derivative 
$\delta(g) := g^{-1}\cdot \dd_{M_n(A)}(g) \in pM_n(\Omega^1_D(A))p$  
leads to 
$$ \nabla^g = \nabla +  \delta(g). $$
More generally, we get 
$$ \nabla_\alpha^g = \nabla +  \delta(g) + \Ad(g^{-1}).\alpha 
= \nabla_{\alpha'}\quad \mbox{ for } \quad 
\alpha'= \delta(g) + \Ad(g^{-1}).\alpha.$$ 
From that we derive in particular that $\nabla_\alpha^g = \nabla_\alpha$ is 
equivalent to 
$\delta(g) = \alpha - \Ad(g^{-1}).\alpha.$ 
\end{remark}

\nin Karl-Hermann Neeb \\ 
Mathematics Department \\ 
Darmstadt University of Technology \\ 
Schlossgartenstrasse 7 \\  
65289 Darmstadt \\ 
Germany \\ 
email: neeb@mathematik.tu-darmstadt.de \\

\end{document}